\documentclass [12 pt]{amsart}
\usepackage{amssymb,latexsym}
\usepackage{xypic}
\theoremstyle{plain}
\newtheorem{Thm}{Theorem}[section]

\newtheorem{Prop}{Proposition}[section]
\newtheorem{Ex}{Example}[section]

\numberwithin{equation}{section}
\theoremstyle{definition}
\newtheorem{Def}{Definition}[section]
\theoremstyle{remark}
\newtheorem{Rem}{\em Remark}

\begin{document}
	
	\title[Grading and Filtrations of Gamma Rings]{Grading and Filtrations of Gamma Rings}
	
	\author{Shadi Shaqaqha and Afnan Dagher}
	
	\address{Department of Mathematics, Yarmouk University, Irbid, Jordan
	}
	\email[Shadi Shaqaqha]{shadi.s@yu.edu.jo}
	\email[Afnan Dagher]{afnand@yu.edu.jo}
	
	

\begin{abstract}
The aim of this paper is to introduce and study graded and filtered gamma rings and gamma modules. We prove that the filtered $\Gamma$-ring (module) is a generalization of the notion of graded ring (module). Also, we construct a graded $\Gamma$-ring from a filtred $\Gamma$-ring. We investigate some properties of graded and fileterd Gamma rings and Gamma modules. Finally we define and study the strongly graded gamma rings.
\end{abstract}
\maketitle
\section{INTRODUCTION}
The concept of $\Gamma$-rings was first introduced by Nobusawa in 1964 (\cite{Nobusawa}).  Later, Barnes generalized the notion of Nobusawa's $\Gamma$-rings and established a new definition of a $\Gamma$-ring (\cite{Barnes}). Nowadays, the notion of $\Gamma$-rings means the $\Gamma$-rings due to Barnes. 
Many mathematicians have been involved in extending the concepts and results of ring theory to the broader framework of the $\Gamma$-ring setting (see e.g. \cite{Anderson, Dumitru, Kyuno, Kyuno2, Kyuno3} and references therein).\\
The graded and filtered rings play a role in various branches of mathematics, and they have become increasingly important. Consequently, the graded
analogues of different concept are widely studied (see \cite{Bahturin1, Cohen, Dade1, Fida, Nastas2, Roozbeh, Shadi}). For example, S. Shaqaqha, in his dissertation, considered gradings and filtrations on colour Lie superalgebras and obtained more Schreier-type formulas in the case of subalgebras of free colour Lie superalgebras (
\cite{Shadi}). So, our  results  are expected to be useful in various applications.\\
In \cite{Fusheng}, a group graded gamma ring is defined to be a $\Gamma$-ring $R$ which is the direct sum of additive subgroups $R_g$, $g\in K$, such that $R_g\Gamma R_h\subseteq R_{gh}$, for all $g, h \in K$, where $K$ is a group defined multiplicatively. In this paper we replace the group grading with a semigroup then we introduce the notion of filterd $\Gamma$-rings. We extend known results in the the case of graded rings to graded $\Gamma$-rings.\\
The article is organized as follows. In Section 2, we recall some definitions and notations that will be used in our work. In Section 3, we introduce the notion of graded gamma rings and we ofer some examples. In Section 4, we obtain some properties of semigroup graded gamma rings. We present ways to obtain new graded gamma rings from old ones. In Section 5, we define filtered gamma rings. We give a recipe to form new graded gamma ring from a filtered gamma ring. In Section 6, we study semigroup graded gamma modules. In Section 7, we obtain a way to produce a ring out of old gamma rings. In Section 8, we study filtered gamma modules and obtain some properties. The strongly graded gamma rings are introduced and studied in Section 9.
\section{Preliminaries}
The notion of gamma rings is a one of the generalizations of rings.
\begin{Def}
Let $\Gamma$ be an additive abelian group. A $\Gamma$-ring (in the sense of Barnes) is an additive abelian group $R$  together with a map
$$R\times \Gamma \times R\rightarrow R : (x,\alpha, y)\mapsto x\alpha y$$
and satisfying the following two identities, for any $x, y, z \in R$ and $\alpha, \beta\in \Gamma$:
\begin{itemize}
\item[(i)] $(x+y)\alpha z= x\alpha z + y\alpha z$, $x(\alpha + \beta) z = x\alpha z+ x\beta z$, $x \alpha (y + z)= x\alpha y + x\alpha z$.
\item[(ii)] $(x\alpha y)\beta z = x\alpha (y \beta z)$.
\end{itemize}
\end{Def}
Clearly every ring $R$ can be regarded as an $R$-ring. Let $G$ be an additive group. We donote by $G_{m,n}$ the set of all $m\times n$ matrices over $G$. Let $M$ be a $\Gamma$-ring. Then $M_{m,n}$ forms a $\Gamma_{n,m}$-ring (see \cite{Kyuno2}).\\
Let $R$ be a $\Gamma$-ring. An element $1\in R$ is called a unity of $R$ if $a\gamma 1 = 1\gamma a = a$ for all $a\in R$ and for some $\gamma\in \Gamma$. In this case, we say $R$ is a $\Gamma$-ring with unity. One can easily show that if a $\Gamma$-ring has a unity with respect to $\gamma\in \Gamma$, then it is unique.\\
A nonempty subset $S$ of a $\Gamma$-ring $R$ is a sub-$\Gamma$-ring of $R$ if $a-b\in S$ and $a\gamma b\in S$ for any $a, b\in S$ and $\gamma\in \Gamma$. A subset $I$ of the $\Gamma$-ring $R$ is a left (right) ideal of $R$ if $I$ is an additive subgroup of $R$ and $R\Gamma I= \{r\alpha a~|~r\in R, \alpha\in \Gamma, a\in I\}$ ($I\Gamma R$) is contained in $I$. If $I$ is
both a left and a right ideal of $I$, then we say that $I$ is an ideal (or two-sided ideal) of $R$. It is clear that the intersection of any number of left (respectively right or two-sided)
ideals of $R$ is also a left (respectively right or two-sided) ideal of $R$.\\
Let $I$ be a (two-sided) ideal of a $\Gamma$-ring $R$. then the factor group $R/I = {x + I:x\in R}$ of all cosets of $I$ is a $\Gamma$-ring with respect to the following multiplication:
\begin{eqnarray*}
(x+I)(\gamma)(y+I)&=& ((x\gamma y)+I
\end{eqnarray*}
for $x, y\in R$ and $\gamma\in \Gamma$.
Let $R_1, R_2, \ldots, R_n$ be $\Gamma$-rings. Then the direct product $\prod_{i=1}^nR_i= R_1\times \cdots\times R_n$ has the structure of a $\Gamma$-ring under the componentwise addition and the following multiplication
$$(r_1, \ldots, r_n)\gamma (s_,, \ldots, s_n)= (r_1\gamma s_1, \ldots, r_n\gamma s_n)$$
for all $r_i, s_i\in R_i~(i=1, \ldots, n)$ and $\gamma\in \Gamma$\\
Suppose that $R_1$ and $R_2$ are $\Gamma$-rings and let $\varphi: \Gamma\rightarrow \Gamma$ be a group isomorphism. Then $f: R_1\rightarrow R_2$ is $\varphi$-homomorphism if $f(x+y)= f(x)+ f(y)$ and $f(x\gamma y)= f(x)\varphi(\gamma) f(y)$ for all $x, y\in R_1$ and $\gamma\in \Gamma$. In particular, if $\varphi$ is the identity map, then $f$ is homomorphism.
\section{Basic Definition and Examples}
\begin{Def}\label{def:graded} Suppose that $G$ is an abelian semigroup written multiplicatively. Let $R$ be a $\Gamma-$ring. Then $R$ is called a graded $\Gamma-$ring of type $G$ if there exist  additive subgroups $R_g$ of $R$ such that $R=\oplus_{g\in G} R_g$ and $R_g \Gamma R_h\subseteq R_{gh}$ for all $g, h\in G$. A nonzero element $r\in R$ is called homogeneous if there is $g\in G$ such that $r\in R_g$, and in this case we say the degree of $r$ is $g$ and we write $d(m)=g$.
\end{Def}
So, any nonzero element $r\in R$, where $R$ is a graded $\Gamma$-ring of type $G$, has a unique representation  as a sum of a finite set  of homogenous elements $r=r_{g_1}+\cdots+ r_{g_k}$.
In this case, we say that $r_{g_1}, r_{g_2}, \ldots$, and $r_{g_k}$ are the homogeneous components of $r$. The set $G_R= \{g\in G~|~R_g\neq\{0\}\}$ is called the support of $R$.\\
Throughout the paper we denote by $G$ an abelian semigroup (written multiplicatively) unless otherwise stated.\\
A two-sided ideal $I$ of a graded $\Gamma$-ring $R$ of type $G$ is called a graded $\Gamma$-ideal (or homogeneous $\Gamma$-ideal) if $I=\bigoplus_{g\in G}I_g$, where $I_g= I\cap R_g$. In other words, $I$ is a graded $\Gamma$-ideal if and only if , for each $x\in I$, its homogeneous components belong also to $I$.
\begin{Ex}\label{exa:graded}
\begin{itemize}
\item[(i)] Let $G$ be a semigroup with identity $e$ (or a monoid). Any $\Gamma$-ring $R$ may be considered as a graded $\Gamma$-ring of type $G$ by putting $R_e= R$, and $R_g=\{0\}$ for any $g\neq e$ in $G$ (i.e. $M$ has the trivial subgroup support). Such a graded $\Gamma$-ring is called trivial.
\item[(ii)] If $I$ is a graded $\Gamma$-ideal of a graded $\Gamma$-ring $R$ of type $G$, then $R/I$ is also a graded $\Gamma$-ring and has the following decomposition:
    $$R/I=\bigoplus_{g\in G}\left(R_g+I\right)/I\cong \bigoplus_{g\in G}\left(R_g/\left(R_g\cap I_g\right)\right).$$
(The last congruence follows from the second isomorphism theorem for groups).
\item[(iii)] Given a semigroup $G$ and a $\Gamma$-ring $R$. Set $RG$ to be the set of all linear combinations
$$r=\sum_{g\in G}r_gg,$$
where $r_g\in R$ and where only finitely many of the $r_g's$ are nonzero. Next, we define the addition and the multiplication on $RG$ as follows: for $\alpha = \sum_{g\in G}a_gg, \beta=\sum_{g\in G}b_gg\in RG$ and $\gamma\in \Gamma$, we have
\begin{eqnarray*}
\alpha +\beta&=& \sum_{g\in G}(a_g+b_g)g,\\
\alpha\gamma\beta&=&\sum_{g,h\in G}(a_g\gamma b_h)gh.
\end{eqnarray*}
Then $RG$ is a graded $\Gamma$-ring of type $G$. Such graded $\Gamma$-ring is called semigroup $\Gamma$-ring (and if $G$ is group, then it is called group Gamma ring). 
\item[(iv)] Let $R$ be a $\Gamma$-ring. Then $R^0=R$ is a $\Gamma$-ring too under the same addition on $R$, but the multiplication is defined by $x\circ \gamma\circ y = y \gamma x$ for $x, y\in R^0$, and $\gamma\in \Gamma$. Furthermore, if $R$ is a graded $\Gamma$-ring of type $G$, where $G$ is a group, then $R^0$ is a graded $\Gamma$-ring of type $G$ too by setting $(R^0)_g = R_{g^{-1}}$, $g\in G$. Indeed for $x\in (R^0)_g$, $y\in (R^0)_h$, and $\alpha\in \Gamma$, we have $x\circ \alpha \circ y= y\alpha x \in R_{h^{-1}g^{-1}}=R_{(gh)^{-1}}=(R^0)_{gh}$.
\item[(v)] Let $R=\bigoplus_{g\in G}R_g$ and $S=\bigoplus_{g\in G}S_g$ be graded $\Gamma$-rings of type $G$. Then the $\Gamma$-ring $R\times S$ may be made into a graded $\Gamma$-ring by setting:
    $$(R\times S)_g = R_g\times S_g~\forall g\in G.$$
More generally, if $(R_i){i\in I}$, where $I$ is finite, is a family of graded $\Gamma$-rings of type $G$, then $R=\prod_{i\in I}R_i$ is a graded $\Gamma$-ring of type $G$ by setting
$$\left(\prod_{i\in I}R_i\right)_g = \prod_{i\in I}(R_i)_{g}.$$
\end{itemize}
\end{Ex}
\section{Some Properties of Graded Gamma Rings}
The following two theorems give us recipes for forming new graded gamma ring from old ones.
\begin{Thm}\label{GGR1}
Let $R$ be a graded $\Gamma$-ring of type $G$ where $G$ is a semigroup. If $\varphi: G\rightarrow H$ is an onto homomorphism (epimorphism) of semigroups, then the $\Gamma$-ring $S= R$ with gradation
$$S_h=\bigoplus_{g\in \varphi^{-1}(h)}R_g$$
for all $h\in H$ is a graded $\Gamma$-ring of type $H$ (such graded $\Gamma$-ring is denoted by $R_{(H)}$).
\end{Thm}
{\it Proof.~} Let $a, b\in H$, $x\in S_a$ and $y\in S_b$. Then $x=x_{a_1}+\cdots + x_{a_r}$ and $y=y_{b_1}+\cdots+y_{b_s}$, where $\varphi(a_i)= a$ and $\varphi(b_j)=b$ for all $i,j$. Now, for $\gamma\in \Gamma$, we obtain
$$x\gamma y= \sum_{1\leq i\leq r, 1\leq j\leq m}x_{a_i}\gamma y_{b_j}.$$
Also, for $1\leq i\leq r$ and $1\leq j\leq m$, we have $\varphi(a_i b_j)=\varphi(a_i)\varphi(b_j)=ab$. Thus $x\gamma y\in S_{ab}$. \hfill $\Box$

Also, we have the following result. We omit the proof since the proof is straightforward.
\begin{Thm}\label{GGR1A}
Let $R$ be a graded $\Gamma$-ring of type $G$, and let $H$ be a subsemigroup of $G$. Then
 \begin{itemize}
 \item[(i)] $R'= \bigoplus_{h\in H}R_h$ is a graded $\Gamma$-ring of type $H$ (such graded is denoted by $R^{(H)}$). In particular, if $e$ is identity of $G$, then $R_e$ corresponds to the trivial subsemigroup of $G$.
 \item[(ii)] If $H$ is a normal subsemigroup of $G$, then
 $$R=\bigoplus_{gN\in G/N}R_{gN}$$
 is a graded $\Gamma$-ring of type $G/N$, where
     $$R_{gN}= \bigoplus_{n\in N}R_{gn}.$$
 \end{itemize}
\end{Thm}
The following theorem extends a well known result to group graded rings.
\begin{Thm}\label{GGR2}
Let $G$ be a monoid with the identity element $e$, and let $R$ be a graded $\Gamma$-ring of type $G$. Then $R_e$ is a sub-$\Gamma$-ring of $R$. Also if $1$ is the unity element with respect to $\gamma_0\in \Gamma$, then $1\in R_e$.
\end{Thm}
{\it Proof.~} For $a, b\in R_e$, we have $a-b\in R_e$ since $M_e$ is an additive subgroup of $R$, and also $a\gamma b\in R_{ee}=R_e$ for any $\gamma\in \Gamma$. It follows that $R_e$ is a sub-$\Gamma$-ring of $R$. Next, as $R=\bigoplus_{g\in G}R_g$, we may assume $1= r_{g_1}+r_{g_2}+\cdots +r_{g_n}$. Pick $\tau \in G$, and a nonzero element $\lambda_\tau \in R_\tau$, then
$$1\gamma_0\lambda_\tau= \lambda_\tau= r_{g_1}\gamma_0\lambda_\tau+r_{g_2}\gamma_0\lambda_\tau+\cdots+ r_{g_m}\gamma_0\lambda_\tau.$$
and so, for all $g_i\in G$ with $g_i\neq e$ we have $r_{g_i}\gamma_0\lambda_\tau = 0$. Hence, $1\gamma_0\lambda_\tau= r_e\gamma_0\lambda_\tau$. Therefore, by the uniqueness of the unity, we have $1=r_e\in R_e$. \hfill$\Box$

Let $R$ be a $\Gamma$-ring with a unity $1$ corresponding to $\gamma_0\in \Gamma$. An element $r\in R$ is called invertible with respect to $1$ and $\gamma_0$ if there is an element (unique) $r^{-1}\in R$ such that $r^{-1}\gamma_0 r = r\gamma_0 r^{-1} = 1$.
\begin{Thm}\label{GGR2A}
Let $G$ be a group with the identity element $e$, and let $R$ be a graded $\Gamma$-ring of type $G$. If $r\in R_g$ is an invertible (corresponding to $\gamma_0\in \Gamma$) homogeneous element, then $r^{-1}$ is homogeneous too of degree $g^{-1}$.
\end{Thm}
{\it Proof.~} Let $s=s_{g_1}+\cdots +s_{g_k}$ be the inverse of $r$. Then
$$1= r\gamma_0 s = r\gamma_0 s_{g_1}+\cdots + r\gamma_0 s_{g_k},$$
where $r\gamma_0 s_{g_i}\in R_{gg_i}$ for all $i=1,\ldots, k$. By Theorem \ref{GGR2}, $1$ is homogeneous of degree $e$, and also, the decomposition is unique, so that $r\gamma_0 s_{g_i}= 0$ for each $g_i \neq g^{-1}$. Since $r$ is invertible, $s_{g^{-1}}\neq 0$, and so $s= s_{g^{-1}}\in R_{g^{-1}}$ as desired. \hfill $\Box$
\section{Filtered Gamma Rings}
The following definition produces the notion of filtered $\Gamma$-rings. It is a generalization of the notion of filtered rings \cite{Nastas}.
\begin{Def}\label{def:filt}
Let $R$ be a $\Gamma$-ring and let $R^0\subseteq R^1\subseteq R^2\subseteq\cdots$ be a chain of additive subgroups of $R$. This chain is called an (ascending) {\em filtration} of $R$ if
$$\bigcup_{r\geq 0}R^m = R ~\mathrm{and} ~ R^i\Gamma R^j\subseteq R^{i+j}~\forall i, j\geq 0.$$
In this case we say that $R$ is a filtered $R_\Gamma$-ring.
\end{Def}
Let $R$ and $S$ be filtered $\Gamma$-rings. A homomorphism $\varphi: R\rightarrow S$ is called a filtered homomorphism if $\varphi(R_i)\subseteq S_i$ for $i=0, 1, \ldots$.
\begin{Ex}\label{FGR1}
\begin{itemize}
\item[(i)] Any $\Gamma$-ring $R$ can be made into a filtered $\Gamma$-ring by setting $R^k=R$ for all $k\geq 0$. Such filtration is called trivial.
\item[(ii)] Let $R$ be a $\Gamma-$ ring, and let
$$R[x]=\{p(x)= a_0+a_1x+\cdots+a_nx^n: a_0,a_1,\cdots,a_n\in R\}.$$
Then it is easy to see that $R[x]$ is a $\Gamma$-ring under the ordinary addition and the multiplication defined as follows: if  $p(x)=a_0+a_1x+\cdots+a_nx^n, q(x)=b_0+b_1x+\cdots+b_mx^m\in R[x]$ such that $a_n \neq 0\neq b_m$, then 
$$p(x).\alpha .  q(x) = a_0 \alpha b_0+ (a_0 \alpha b_1+b_0 \alpha a_1)x+ \cdots+ a_n\alpha b_m x^{m+n}$$. Moreover, $R[x]$ has a filtration $\bigcup_{k=0}^\infty R^k$, where $R^k$ is spanned by all monomials of degree less than or equal to $
k$.
\item[(iii)] Given a graded $\Gamma$-ring $R=\bigoplus_{i\geq 0}R_i$ of type $\mathbb{Z}$, then there is a corresponding filtration $$\bigcup_{k=0}^\infty R^k,$$
where $R^k=\bigoplus_{j=0}^kR
_j$.
\end{itemize}
\end{Ex}
Conversely, the following result gives us a way to produce a graded $\Gamma$-ring out of a filtered $\Gamma$-ring.
\begin{Thm}\label{GGR3}
Let $R$ be a $\Gamma$-ring with filtration $R^0\subseteq R^1\subseteq R^2 \subseteq \cdots$. Set $R^{-1}=0$, and let $\mathrm{gr}R = \bigoplus_{k=0}^{\infty}\left(R^k/R^{k-1}\right)$. Then define multiplication
$$R^m/R^{m-1}\times \Gamma\times R^n/R^{n-1}\rightarrow R^{m+n}/R^{m+n-1}$$
by $\left(r+R^{m-1}\right)\gamma \left(s+R^{n-1}\right)= (r\gamma s)+R^{m+n-1}$. Extending this multiplication by multilinearity makes $\mathrm{gr}R$ is a graded $\Gamma$-ring.
\end{Thm}
{\it Proof.~}We must show that the map is well defined. Suppose that $x+R^{m-1}=x\ '+R^{m-1}$, $y+R^{n-1}=y\ '+R^{n-1}$, and $\gamma\in \Gamma$. Then $x-x\ '= r_0$ and $y-y\ '=r_1$ for some $r_0\in R^{m-1}$ and $r_1\in R^{n-1}$. Thus
\begin{eqnarray*}
x\gamma y -x\ '\gamma y\ '&=& (x\ '+r_0)\gamma (y\ '+r_1)- x\ '\gamma y\ '\\
&=& x\ '\gamma r_1+ r_0\gamma y\ '+ r_0\gamma r_1
\end{eqnarray*}
is a member of $R^{m+n-1}$. The remainder of the proof is quite straightforward.

\hfill$\Box$
\begin{Rem}\label{GGR4}
If the filtration $R^k$, given in the theorem above, comes from a grading (i.e. $R^k = \bigoplus_{i=0}^kR_i$), then $R^k/R^{k-1}\cong R_k$ as additive groups, and so $\mathrm{gr}R\cong R$ as $\Gamma$-rings.
\end{Rem}
\section{Graded Gamma Modules}
Let $R$ be a $\Gamma$-ring. An abelian group $M$ together with a mapping
$$R\times \Gamma\times M\rightarrow M;~ (r, \gamma, m)\mapsto r\gamma m,$$
such that the following identities are satisfied for all $r_1, r_2, r\in R$, $\gamma_1, \gamma_2, \gamma\in \Gamma$, and $m_1, m_2, m\in M$:
\begin{itemize}
	\item[(i)] $r\gamma(m_1+ m_2)= r\gamma m_1+ r\gamma m_2$,
	\item[(ii)]  $(r_1+ r_2)\gamma m= r_1\gamma m + r_2\gamma m$,
	\item[(iii)] $r(\gamma_1+ \gamma_2)m= r\gamma_1 m + r\gamma_2 m$,
	\item[(iv)] $(r_1\gamma_1)(r_2\gamma_2 m)= (r_1\gamma_1 r_2)\gamma_2 m$
\end{itemize} 
is called a (left) $R_\Gamma$-module. A right $R_\Gamma$-module is defined in analogous manner (\cite{Ameri}). If $R$ and $S$ are $\Gamma$-rings, then we say $M$ is an $(R, S)_\Gamma$-bimodule if it is both left $R_\Gamma$-module and right $S_\Gamma$-module and simultaneously $(r\alpha m)\beta s = r\alpha (m\beta s)$ for all $r\in R$, $m\in M$, $s\in S$, and $\alpha, \beta \in \Gamma$.
A (left) $R_\Gamma$-module $M$ is unitary if there exists an element , say $1$, in $R$ and $\gamma_0\in\Gamma$ such that $1 \gamma_0 m=m~\forall m\in M$. Note that if $M$ is a left $R_\Gamma$-module then it is easy to verify that $0\gamma m = r0m = r\gamma0 = 0$ for any $r\in R, \gamma\in \Gamma$ and $m\in M$.
\begin{Ex}\label{GM2}
\begin{itemize}
\item[(i)] If $R$ is a $\Gamma$-ring and $M$ is an abelian group, then $M$ is an $R_\Gamma$-module by letting $r\gamma m =0~\forall r\in R, \gamma\in \Gamma, m\in M$.
\item[(ii)] Every $\Gamma$-ring $R$ is an $R_\Gamma$-module by
$$R\times \Gamma\times R\rightarrow R;~(r_1, \gamma, r_2)\mapsto r_1\gamma r_2.$$
\item[(iii)] Let $R$ be a graded $\Gamma$-ring of type $G$, where $G$ is a monoid with identity $e$. Then $R$ can be considered as an $\left(R_e, R_e\right)_{\Gamma}$-bimodule.
\end{itemize}
\end{Ex}
Let $R$ be a $\Gamma$-ring. A nonempty subset $M_1$ of $M$ is called a (left) $R_\Gamma$-submodule of $M$ if $M_1$ is a subgroup of $M$ and $R\Gamma M_1\subseteq M_1$. In this case we write $M_1\leq M$. For a $\Gamma$-ring $R$, $\{0\}$ and $R$ are $R_\Gamma$-submodules of $R_\Gamma$-module $R$. Also, every ideal of $R$ is a submodule of $R$.\\
The following definition is a generalization of the notion of a graded module.
\begin{Def}\label{GM3}
Let $R$ be a graded $\Gamma$-ring of type $G$, and $M$ be an $R_\Gamma$-module. Then $M$ is said to be graded (left) $R_\Gamma$-module if there is a family $\{M_g~:~g\in G\}$ of additive subgroups of $M$ such that $M=\bigoplus_{g\in G}M_g$ and $R_g\Gamma M_h\subseteq M_{gh}$ for all $g, h\in G$. A nonzero element $m\in M$ is called homogeneous of degree $g$, and we write $\mathrm{deg}(m)=g$, if there $g\in G$ such that $m\in M_g$. A submodule $N$ of $M$ is a graded submodule if $N=\bigoplus_{g\in G}(N\cap M_g)$ (equivalently, for any $x\in N$ the homogeneous component(s) of $x$ are again in $N$).
\end{Def}
\begin{Ex}
Let $R$ be a graded $\Gamma$-ring of type $G$, $M$ be a graded $R_\Gamma$-module, and $K$ be a submodule of $M$. Then $K'=\bigoplus_{g\in G}K_g$, where $K_g$ is the additive subgroup of $K$ generated by $K\cap (h(M))_g$ where $(h(M))_g$ is the set of all homogeneous elements of degree $g$ (i.e. $K_g$ is the smallest subgroup of $G$ that contains $K\cap (h(M))_g$), is a graded submodule of $M$. One can easily show that it is a maximal submodule of $K$ that is a graded submodule of $M$.
\end{Ex}
It is known that if $R$ is a $\Gamma$ ring, $M$ is a $R_\Gamma$-ring, and $K$ is a subgroup of $M$, then that factor group $M/K$ is an $R_\Gamma$-module under the mapping
$$R\times \Gamma \times M/K\rightharpoonup M/K;~(r, \gamma, m+K)\mapsto (r\gamma m)+ M$$
(See \cite{Ameri}). Using the same arguments in the graded rings we can prove the following result.
\begin{Prop}\label{GM4}
Let $R$ be a graded $\Gamma$-ring of type $G$ and $M$ be a graded $R_\Gamma$-module. If $K$ is a submodule of $M$, then the factor module $M/K$ is a graded $R_\Gamma$-module by setting
$$(M/K)_g =(M_g+K)/K\cong M_g/(K\cap M_g)$$
for any $g\in G$.
\end{Prop}
\section{On Homomorphisms of Gamma Modules}
Suppose that $R$ is a graded $\Gamma$-ring of type $G$ and that $M$ and $K$ are (left) graded $R_{\Gamma}$-module. Let $\varphi: \Gamma\rightarrow \Gamma$ be a group isomorphism. Then $f: M\rightarrow K$ is $\varphi$-homomorphism if $f(x+y)= f(x)+ f(y)$ and $f(r\gamma x)= r\varphi(\gamma) f(x)$ for all $x, y\in M$, $r\in R$, and $\gamma\in \Gamma$. In particular, if $\varphi$ is the identity map, then $f$ is homomorphism. It is called isomorphism if it is one to one and onto.\\
For graded (left) $R_{\Gamma}$-modules $M$ and $K$. A ($\varphi$-)homomorphism $f: M\rightarrow K$ is called homogeneous of degree $h\in G$ if for all $g\in G$, we have $f(M_g)\subseteq K_{hg}$. In particular, for graded $\Gamma$-rings $R$ and $K$. A ($\varphi$-)homomorphism of $\Gamma$-ringe $f: R\rightarrow K$ is called homogeneous of degree $h\in G$ if for all $g\in G$, we have $f(R_g)\subseteq K_{hg}$. A homogeneous $(\varphi-)$homomorphism of $R_\Gamma$-modules of degree $e\in G$ (the identity element of $G$) will be called a degree preserving map. Clearly the set  $\mathrm{Hom}_R(M, K)$ of all homomorphisms of $R_\Gamma$ modules from $M$ into $K$ forms a group under the ordinary addition. Also, for graded $R_{\Gamma}$-modules $M_1, M_2$, and $M_3$ the composition of a homogeneous ($\varphi$-)homomorphism $f: M_1\rightarrow M_2$ of degree $h\in G$ with a ($\varphi$-)homomorphism $g: M_2\rightarrow M_3$ of degree $k\in G$ is ($\varphi$-)homomorphism of degree $hk$.
\begin{Def}\label{G5}
Let $R$ be a $\Gamma$-ring and $M$ be an $R_{\Gamma}$-module. Then we say $M$ is finitely generated if there exist $m_1, \ldots, m_k$ in $M$ such that if $x\in M$, then there exist $r_1, \ldots, r_k\in R$ and $\gamma_1, \ldots, \gamma_k\in \Gamma$ with $x=r_1\gamma_1m_1+\cdots + r_k\gamma_k m_k$. In this case the set $\{m_1, \ldots, m_k\}$ is said a generating set to $M$.
\end{Def}
The following theorem gives us a procedure to produce a ring out of old $\Gamma$-ring.
\begin{Thm}\label{GM6}
Let $M$ and $K$ be graded $R_{\Gamma}$-modules of type $G$. If $M$ is a finitely generated, then
$$\mathrm{Hom}_R(M, K)=\bigoplus_{g\in G}\left(\mathrm{Hom}(M, K)\right)_g$$
is a graded abelian group, where $\mathrm{Hom}_R(M, K)_g$  is the subgroup of $\mathrm{Hom}_R(M, K)$ that consisting of all homomorphisms of from $M$ to $K$ of degree $g$. Moreover if $M=K$, then $\mathrm{Hom}(M, M)$ is a graded ring of type $G$.
\end{Thm}
{\it Proof.~} Suppose that $f\in \mathrm{Hom}_R(M, K)$. Let $g\in G$. Define a map $$f_g:M \rightarrow K$$
as follows:
for $m= m_{g_1} + m_{g_2}+ \cdots + m_{g_k}\in M$ ($g_1, g_2, \ldots, g_k\in G$), set
$$f_g(m)= \left(f(m_{g_1g^{-1}})\right)_{g_1}+\cdots+ \left(f(m_{g_kg^{-1}})\right)_{g_k}.$$
Then $f_g\in \mathrm{Hom}_R(M, K)$. Indeed if $x=\sum_{h\in G}x_h$, $y=\sum_{h\in G}y_h\in M$, then
\begin{eqnarray*}
f_g(x+y)&=&\sum_{h\in G}\left(f(x+y)_{hg^{-1}}\right)_h\\
        &=&\sum_{h\in G}\left((f(x_{hg^{-1}}+y_{hg^{-1}})\right)_h\\
        &=&\sum_{h\in G}\left((f(x_{hg^{-1}})+ f(y_{hg^{-1}})\right)_h\\
        &=&\sum_{h\in G}\left(\left(f(x_{hg^{-1}}\right)_h+\left((f(y_{hg^{-1}}\right)_h\right)\\
        &=&f_g(x)+f_g(y).
\end{eqnarray*}
Also, for $r\in R$ and $\gamma\in \Gamma$, we have
\begin{eqnarray*}
f_g(r\gamma x)&=& \sum_{h\in G}\left(f((r\gamma x)_{hg^{-1}})\right)_h\\
              &=& \sum_{h\in G}\left(r\gamma f(x_{hg^{-1}})\right)_h\\
              &=&r\gamma \sum_{h\in G}\left((f(x_{hg^{-1}})\right)_h\\
              &=&r\gamma f_g(x).
\end{eqnarray*}
In addition, $f_g\in \left(\mathrm{Hom}_R(M, K)\right)_g$. Indeed for any $h\in G$ and $m\in M_h$, we have $f_g(m)= (f(m))_{hg}\in K_{hg}$.\\
Let us assume that $\{m_1, \ldots, m_s\}$ is a generating set to $M$. For $x\in M$, there exist $r_1, \ldots, r_s\in R$ and $\gamma_1, \ldots, \gamma_s\in \Gamma$ with $x=r_1\gamma_1m_1+\cdots + r_s\gamma_s m_s$. Thus $f(x)= r_1\gamma_1f(m_1)+\cdots+ r_s\gamma_sf(m_s)$. On the other hand for each $t=1, 2, \ldots, r$, $f_g(m_t)= 0$ for all but a finite number of $g\in G$ and
$$\sum_{g\in G}f_g(m_t)= \sum_{g\in G}(f(m_t))_{hg}=f(m_t).$$
It follows again that $f_g(x)=0$ for all but finitely number of $g\in G$, and $f= \sum_{g\in G}f_g$. We have so far shown that $\mathrm{Hom}_R(M, K)=\sum_{g\in G}(\mathrm{Hom}_R(M, K))_g$. On the other hand if $f\in \mathrm{Hom}_R(M, K)_{h_1}\cap  \mathrm{Hom}_R(M, K)_{h_2}$ where $h_1, h_2\in G$ and $h_1\neq h_2$, then for any $g\in G$ we have $f(M_g)\subseteq K_{h_1g}\cap K_{h_2g}=\{0\}$. Hence $f=0$. This shows that $\mathrm{Hom}_R(M, K)= \bigoplus_{g\in G}(\mathrm{Hom}_R(M, K))_g$. In particular, if $M=K$, we obtain $\mathrm{Hom}_R(M, M)= \bigoplus_{g\in G}(\mathrm{Hom}_R(M, M))_g$. For $\rho\in \mathrm{Hom}_R(M, M))_{h_1}$ and $\sigma\in \mathrm{Hom}_R(M, M)_{h_2}$ we have $\rho\sigma\in \mathrm{Hom}_R(M, M)_{h_1h_2}$. It is clear now that $\mathrm{Hom}_R(M, M)$ is a graded ring of type $G$.\hfill $\Box$
\section{Filtered $\Gamma$-Modules}
Let $R$ be a filtered $\Gamma$-ring and $M$ be a (left) $R_{\Gamma}$-module. Then $M$ is called a filtered $R_\Gamma$ module if there is an (ascending) chain
$$M^0\subseteq M^1\subseteq M^2 \subseteq \cdots$$
of additive subgroups of $M$ such that $\bigcup_{k\geq 0}M^k= M$ and $R^i\Gamma M^j\subseteq M^{i+j}$ for any $i, j$.
\begin{Ex}
\begin{itemize}
\item[(i)] It is clear that if $R$ is a filtered $\Gamma$-ring, then $R$ is a filtered $R_{\Gamma}$-module.
 \item[(ii)] If the filtration of a graded $\Gamma$-ring $R$ is trivial and $M$ is $R_{\Gamma}$-module, then any ascending chain of submodules $K^s$, where
 $\bigcup_{s\geq 0}K^s=M$, of $M$ defines a filtration for $M$.
 \item[(iii)] Let $R$ be a $\Gamma$-ring and $I\subseteq R$ be an ideal of $R$. Then $R$ has a descending filtration $\bigcup_{i=0}^{\infty}R^k$, where
 $$R^0=R, R^1= I, R^2=I\Gamma I,~ \mathrm{and}~ R^k= R^{k-1}\Gamma I~\mathrm{for}~ k\geq 3.$$
 Such filtration is called $I$-adic filtration. Let $M$ be $R_{\Gamma}$-module. The corresponding descending filtration of $M$
 $$M\supseteq R^1\Gamma M\supseteq  R^2\Gamma M\supseteq \cdots.$$
\end{itemize}
\end{Ex}
Let $R$ be a graded $\Gamma$-ring and let $M=\bigoplus_{i=0}^\infty M_i$ be a graded $R_\Gamma$-module. Consider the corresponding filtration of $R$ as in Example \ref{FGR1}(iii), then there is 
a corresponding filtration on $M$
$$M^0\subseteq M^1\subseteq M^2 \subseteq \cdots$$
where $M^k=\bigoplus_{j=0}^k M_j$. Conversely, one can prove the following result.
\begin{Thm}\label{FGM2}
If $M$ is a filtered $R_\Gamma$ module, then 
$$\mathrm{gr}(M)=\bigoplus_{i=0}^\infty M_i/M_{i-1}$$
is a graded $(\mathrm{gr}R)_\Gamma$ module where $\mathrm{gr}R=\bigoplus_{i=0}^\infty R_i/R_{i-1}$.
\end{Thm}
Let $R$ be a filtered $\Gamma$-ring and $M=\bigcup_{k\geq 0}M^k$ be a filtered $R_{\Gamma}$-module. It is trivial that $\bigcup_{k\geq}M^k$ is an $R_{\Gamma}$-submodule of $M$. The following theorem answers the question: Is $\bigcap_{k\geq}M^k$ an $R_{\Gamma}$-submodule?
\begin{Thm}\label{FGM1}
Let $R$ be a filtered $\Gamma$-ring and $M$ be a filtered $R_{\Gamma}$-module by $M^k$. Then $\bigcap M^k$ is an $R_{\Gamma}$-submodule of $M$.
\end{Thm}
{\it Proof.~} The intersection of any arbitrary collection of subgroups of a group is subgroup, so $\bigcap M^k$ is a subgroup of $M$. For $\lambda\in R$, $\gamma\in \Gamma$, and $x\in \bigcap M^k$. There exists $p\in \mathbb{N}_0$ such that $\lambda \in R^p$. As $x\in M^{k-p}$ for all $k$, we have $\lambda\gamma x\in M^k$ for all $k$. That is $\lambda \gamma x\in \bigcap M^k$. \hfill $\Box$
\begin{Rem}\label{GM7}
	One can focus on studying graded $\Gamma$-bimodules. Let $R$ and $S$ be graded $\Gamma$-rings of type $G$ and let $M$ be an abelian group. Then $M=\bigoplus_{g\in G}M_g$ is a graded $(R, S)_{\Gamma}$-bimodule of type $G$ if the following conditions are satisfied:
	\begin{itemize}
		\item[(i)] $M$ is an $(R, S)_{\Gamma}$-bimodule,
		\item[(ii)] $M$ is a graded left $R_{\Gamma}$-module,
		\item[(iii)] $M$ is a graded right $S_{\Gamma}$-module.
	\end{itemize}
	Therefore $N_g\Gamma M_h\Gamma S_k\subseteq M_{ghk}$ for $g, h, k\in G$.
\end{Rem}
\section{Strongly Graded Gamma Rings}
Let $M$ be a graded $\Gamma$-ring, where $G$ is a monoid, with identity $1$. According to Theorem \ref{GGR2}, we have $1\in M_{e}$. Thus $M_g\Gamma M_e=M_g $ and $M_e\Gamma M_g= M_g$ for all $g\in G$. If these equalities hold for any arbitrary two elements $g, h, G$, we get a strongly graded $\Gamma$-ring. In other words we have the following definition.
\begin{Def}\label{SGGR1}
Let $G$ be a semigroup. A graded $\Gamma$-ring $M$ of type $G$ is called a strongly graded $\Gamma$-ring if $M_g\Gamma M_h= M_{gh}$ for all $g, h\in G$.
\end{Def}
\begin{Ex}\label{SGGR1A}
\begin{itemize}
\item[(i)] Let $G$ be a group and $M$ be a $\Gamma$-ring. The group $\Gamma$-ring $MG$ is a strongly graded $\Gamma$-ring of type $G$.
\item[(ii)] If $M$ is a strongly graded $\Gamma$-ring, then so is $M/I$ for every graded ideal $I$ of $M$.
\item[(iii)] If $M$ is a strongly graded $\Gamma$-ring of type $G$ and $\phi :G \rightarrow H$ is an epimorphism, then $M_{(H)}$ is strongly graded.
\item[(iv)] If $M$ is a strongly graded $\Gamma$-ring of type $G$, and $H$ is a sub(semi-)group of $G$, then $M^{(H)}$ is strongly graded.
\item[(v)] Let $(M_i)_{i\in I}$ be a family of graded $\Gamma$ rings of type $G$ such that $I$ is finite. Then $M=\prod_{i\in I}M_i$ is strongly graded if and only if $M_i$ is strongly graded for each $i\in I$.
\end{itemize}
\end{Ex}
Let $M$ be a $\Gamma$-ring with unity $1_{\gamma_0}$, and let $M^*$ denotes the set of invertible elements in $M$ with respect to $\gamma_0$. Define multiplication on $M^*$ as follows: $mn = m\gamma_0 n$ for $m, n\in M$. Then the set $M^*$ forms a group. For if $m, n\in M^*$, then $mnn^{-1}m^{-1}= (m\gamma_0 n )\gamma_0 (n^{-1}\gamma_0 n^{-1})= 1_{\gamma_0}$. Also, the associativity follows from one of the properties of $\Gamma$-rings. The support of invertible homogeneous elements of $M$ is defined by
$$G_M^*=\{g\in G~|~M_g^*\neq \phi\}$$
where $M_g^*= M_g\cap M^*$. Using Theorem \ref{GGR2A} it is easy to see that $G_M^*$ is a group. Also $G_M^*\subseteq G_M$ where $G_M= \{g\in G~|~M_g\neq \{0\}\}$ is the support of $M$.
\begin{Def}\label{SGGR2}
A graded $\Gamma$-ring $M=\bigoplus_{g\in G}M_g$ with unity $1_{\gamma_0}$, where $G$ is a group, is called crossed product if there is an invertible element in every homogeneous element of $M$, that is $G_M^*= G$.
\end{Def}
\begin{Thm}\label{SGGR3}
Let $M=\bigoplus_{g\in G}M_g$ be a graded $\Gamma$-ring of type $G$, where $G$ is a group, with unity $1_{\gamma_0}$. Then
\begin{itemize}
\item[(i)] $M$ is strongly graded if and only if $1_{\gamma_0}\in M_g\Gamma M_{g^{-1}}$ for all $g\in G$,
\item[(ii)] if $M$ is strongly graded, then the support of $M$ is $G$,
\item[(iii)] any crossed product $\Gamma$-ring is strongly graded,
\item[iv)] if $N$ is a graded $\Gamma$-ring and $\varphi:M\rightarrow N$ is a degree preserving homomorphism of $\Gamma$-rings, and $M$ is strongly graded $\Gamma$-ring, then $N$ is strongly graded too.
\end{itemize}
\end{Thm}
{\it Proof.~} \begin{itemize}
\item[(i)] Suppose that $M$ is a strongly graded $\Gamma$-ring. Then $M_e=M_g\Gamma M_{g^{-1}}$ for any $g\in G$. Now the result follows from Theorem \ref{GGR2}. Conversely, let $x\in M_e$ and $g\in G$. Then $x=1_{\gamma_0}\gamma_0x\in \left(M_g\Gamma M_{g^{-1}}\right)\gamma_0M_e= M_g\Gamma\left(M_{g^{-1}}\gamma_0M_e\right)\subseteq M_g\Gamma M_{g^{-1}}$. Hence $M_e = M_g\Gamma M_{g^{-1}}$ for all $g\in G$. Next for $g, h\in G$ we have
\begin{eqnarray*}
M_{gh} &\subseteq & M_e\Gamma M_{gh}\\
       &=&\left(M_g\Gamma M_{g^{-1}}\right)\Gamma M_{gh}\\
       &=&M_g\Gamma \left(M_{g^{-1}}\Gamma M_{gh}\right)\\
       &\subseteq & M_g\Gamma M_h.
\end{eqnarray*}
This proves that $M_{gh}= M_g\Gamma M_h$ for all $g, h\in G$. Therefore $M$ is strongly graded.
\item[(ii)] Let $g\in G$. Then using (i), we have $1_{\gamma_0}\in M_g\Gamma M_{g^{-1}}$. Thus $M_g\neq \{0\}$.
\item[(iii)] Suppose that $M$ is a crossed product $\Gamma$-ring. Let $g\in G$. Then there exists an invertible element (homogeneous) $m\in M_g$. According to Theorem \ref{GGR2A} we have $m^{-1}\in M_{g^{-1}}$. So $1_{\gamma_0}=m\gamma_0 m^{-1}\in M_g\Gamma M_{g^{-1}}$. Now the result follows from (i).
\item[(iv)] $M$ is strongly graded, so for all $g\in G$ we have $1_{\gamma_0}\in M_g\Gamma M_{g^{-1}}$. It follows $\varphi(1_{\gamma_0})= 1_{\gamma_0}\in \varphi(M_g)\Gamma \varphi(M_{g^{-1}})\subseteq N_g\Gamma N_{g^{-1}}$ for all $g\in G$.
\end{itemize}
\hfill $\Box$
\begin{Rem}
Let $R$ be a graded $\Gamma$-ring, and $M$ be a graded $R_{\Gamma}$-module. Then $M$ is strongly graded $R_{\Gamma}$-module if $R_g\Gamma M_h=M_{gh}$ for all $g, h\in G$. One can focus on considering such graded $\Gamma$-rings. Also, one of the future works will focus on studying tensor product of graded $\Gamma$-modules.
\end{Rem}

\end{document}